\documentclass[11pt]{amsart}

\usepackage{xcolor}
\definecolor{InternalLinkColor}{HTML}{1F4E79}
\definecolor{CitationLinkColor}{HTML}{387043}
\definecolor{URLLinkColor}{HTML}{9C3D10}

\usepackage{hyperref}
\hypersetup{
  colorlinks=true,
  linkcolor=InternalLinkColor,
  citecolor=CitationLinkColor,
  urlcolor=URLLinkColor
}
\usepackage[T1]{fontenc}
\usepackage{lmodern,microtype}
\usepackage{comment}
\usepackage{amsmath,amssymb,mathtools}
\usepackage[margin=1in]{geometry}
\newtheorem{theorem}{Theorem}[section]

\newtheorem{lemma}[theorem]{Lemma}

\newtheorem{corollary}[theorem]{Corollary}
\newtheorem{conjecture}[theorem]{Conjecture}

\theoremstyle{definition}
\newtheorem{definition}[theorem]{Definition}

\newcommand{\T}{\mathbb T}
\newcommand{\R}{\mathbb R}
\newcommand{\Z}{\mathbb Z}
\newcommand{\Q}{\mathbb Q}
\newcommand{\cS}{\mathcal S}
\newcommand{\cV}{\mathcal V}
\newcommand{\cU}{\mathcal U}

\newcommand{\Prog}{\operatorname{Prog}}

\newcommand{\acc}{\operatorname{acc}}
\newcommand{\ml}{\operatorname{ML}}
\newcommand{\mat}[1]{\mathbf{#1}}
\newcommand{\vecb}[1]{\boldsymbol{#1}}

\title{Arithmetic structure of Lonely Runner spectra}
\author{Tanupat Trakulthongchai}
\date{\today}
\subjclass[2020]{Primary 11J13; Secondary 11J06, 11B25, 52C07}
\keywords{Lonely Runner Conjecture, view-obstruction problem, rational piecewise-affine functions}
\begin{document}

\begin{abstract}
For a fixed integer $n\geq1$, let $D(T)$ be the $L^\infty$-distance from a subtorus $T\subseteq(\R/\Z)^n$
to $(1/2,\ldots,1/2)$. For a subtorus $U\subseteq(\R/\Z)^n$ and an integer $1\leq k\leq n$, its $k$th-order \emph{relative Lonely Runner spectrum} $\mathcal{S}_k(U)$ is defined as the set of $D(T)$ as $T$ ranges over all $k$-dimensional proper subtori of $U$. We show that when $\dim U=k+1$, the spectrum $\mathcal{S}_k(U)-D(U)$ is the union of the reciprocals of finitely many arithmetic progressions up to finite symmetric difference, extending a result of Jain and Kravitz for $k=1$. In fact, we prove a more general statement concerning the approximation of the minima of rational piecewise-affine functions.
\end{abstract}

\maketitle

\section{Introduction}
\subsection{The Lonely Runner Conjecture and spectra}
At the start line of a unit-length circular track, $n$ runners and a referee are standing. When the referee signals, the runners begin to run indefinitely with constant but distinct non-zero integer speeds $v_1,\ldots,v_n$. How far can the stationary referee be from all of the runners simultaneously? This length is called the \emph{maximum loneliness} $\ml(v_1,\ldots,v_n)$. Wills' \cite{wills_zwei_1967} conjecture that it is always at least $1/(n+1)$ has come to be known as the Lonely Runner Conjecture, due to the above popular interpretation attributed to Goddyn \cite{bienia_flows_1998}. At the time of writing, it has been verified for $n\leq14$ \cite{betke1972,cusick1984view,bohman2001six,barajas_lonely_2008,rosenfeld2026lonely,trakulthongchai2026nine,sungkawichai2026eleven,allikvere2026fourteen}. A more complete account of the conjecture was given in \cite{perarnau_lonely_2025}, although it was published before the recent advances that resolve the conjecture for $7\leq n\leq14$.

Cusick \cite{cusick_view-obstruction_1973} reinterpreted this conjecture as a view-obstruction problem, which we restate in the language of Giri and Kravitz \cite{giri_structure_2026}. Let $\T=\R/\Z$ and let $\|\cdot\|_p$ denote the standard $p$-norm for $p\in[1,\infty]$ (when this norm is applied on a torus we choose a representative that makes the norm smallest). For $\vecb{x}\in\T^n$, define
\[
 D(\vecb{x})=\|\vecb{x}-(1/2,\ldots,1/2)\|_\infty\] and $D(X)=\min_{\vecb{x}\in X}D(\vecb{x})$ for non-empty and closed $X\subseteq\T^n$. It is easily seen that
\[\ml(v_1,\ldots,v_n)=1/2-D(\langle(v_1,\ldots,v_n)\rangle_\R),\] so the Lonely Runner Conjecture is equivalent to the statement that $D(\langle\vecb{v}\rangle_\R)\leq\frac{1}{2}-\frac{1}{n+1}$ for every $\vecb{v}\in\Z_{>0}^n$.

Other than asking for the maximum of $D(\langle\vecb{v}\rangle_\R)$, one could ask for all possible values that $D(\langle\vecb{v}\rangle_\R)$ can attain. The first substantial work in this direction is due to Kravitz \cite{kravitz_barely_2021}, who anticipated the following dichotomy: for each speed tuple $\vecb{v}\in\Z_{>0}^n$, the runners are either ``very lonely'' with $D(\langle\vecb{v}\rangle_\R)\leq\frac{1}{2}-\frac{1}{n}$, or ``barely lonely'' with $D(\langle\vecb{v}\rangle_\R)=\frac{1}{2}-\frac{s}{ns+1}$ for $s\in\Z_{>0}$. If this was true, the quantity $\frac12-\frac{1}{n+1}$ in the Lonely Runner Conjecture could be explained as the first element of a structured arithmetic sequence. This statement was proven for $n=3$ and partially for $n=4,6$ in \cite{kravitz_barely_2021}. However, Fan and Sun \cite{fan_amending_2026} refuted it for $n=4$, by showing that one can also have 
\begin{equation}
\label{eq:fansun}
    D(\langle\vecb{v}\rangle_\R)=\frac{1}{2}-\frac{1}{4}+\frac{1}{16s+60}
\end{equation}
for each $s\in\Z_{\geq0}$. Nevertheless, this sequence still converges to the same limit of $\frac{1}{2}-\frac{1}{n}$ as the original sequence. Motivated by this observation, Giri and Kravitz  \cite{giri_structure_2026} studied the accumulation of the \emph{Lonely Runner spectra}.

A subtorus $U$ is a connected and closed subgroup of $\T^n$, and is said to be \emph{proper} if it is not contained in the union of the coordinate hyperplanes, or equivalently if $D(U)<1/2$. An example of a subtorus is $H_{\vecb{\lambda}}$ where $\vecb{\lambda}\in\Z^n$ is primitive (corresponding to connectedness), which is defined as \[
H_{\vecb{\lambda}}=\{\vecb{x}\in\T^n:\vecb{\lambda}\cdot\vecb{x}\equiv0\pmod{1}\}.
\] For a proper subtorus $U$ and $1\leq k\leq\dim U$, its $k$th-order \emph{Lonely Runner spectrum} relative to $U$ is given by
\[
 \cS_k(U)=\{D(T):T\subseteq U\text{ is a }k\text{-dimensional
 proper subtorus}\}.
\]
The \emph{classical} Lonely Runner spectra are $\cS_k(\T^n)$. The most familiar is $\cS_{1}(\T^n)$, the set of values of $D(\langle\vecb{v}\rangle_\R)$. We write $\acc(S)$ for the set of accumulation points of $S$. 

In \cite{giri_structure_2026}, $\cS_k(\T^n)$ was shown to accumulate only from above, and the inclusions
\begin{equation}
    \label{eq:specincl}
    \cS_1(\T^{n-1})\subseteq\acc(\cS_1(\T^n))\subseteq\cS_2(\T^{n}),
\end{equation}
and $\acc(\cS_k(\T^n))\subseteq\cS_{k+1}(\T^n)$ were established. Several conjectures were then posed (almost all natural generalizations of the previous results are open questions), but the following is most relevant to us.
\begin{conjecture}[{\cite[Problem 1.6]{giri_structure_2026}}]
\label{conj:spectra}
    For any $1\leq k\leq n$, we have
    $\cS_k(\T^n)=\cS_{k+1}(\T^{n+1})$.
\end{conjecture}
In other words, it was conjectured that the spectra $\cS_k(\T^n)$ depend only on the codimension $n-k$. Conjecture \ref{conj:spectra} is appealing because, by inclusion \eqref{eq:specincl}, it implies the accumulation conjecture \cite[Conjecture 1.5]{giri_structure_2026} that $\acc(\cS_1(\T^n))=\cS_1(\T^{n-1})$, which depicts a nice relationship between different dimensions of the Lonely Runner problem. Together with the Lonely Runner Conjecture in dimension $n-1$, this would imply, for example, that any sequence of distinct values of $D(\langle\vecb{v}\rangle_\R)$ larger than $\frac12-\frac{1}{n}$ converges to $\frac12-\frac{1}{n}$. The inclusion $\cS_k(\T^n)\subseteq\cS_{k+1}(\T^{n+1})$ is obvious; the difficulty lies in showing the opposite inclusion. We can do this for $n-k=1$ because\footnote{See \cite[Section 8]{giri_structure_2026} for the proof, or simply note for $\vecb{\lambda}$ primitive, $D(H_{\vecb{\lambda}})=\frac{1}{2\|\vecb{\lambda}\|_1}$ if $\|\vecb{\lambda}\|_1$ is odd, and $D(H_{\vecb{\lambda}})=0$ if $\|\vecb{\lambda}\|_1$ is even.} \[\cS_{n-1}(\T^n)=\{0\}\cup\{1/(4s+2):s\in\Z_{>0}\},\] but it remains wide open for $n-k>1$.

\subsection{The structure theorems}
The results in \cite{giri_structure_2026} were stated in the form of $\cS_k(\T^n)$, but Jain and Kravitz \cite{jain_relative_2026} noticed that most of them still hold when $\T^n$ is replaced with an arbitrary $U$ with dimension $n$. Therefore, they proposed studying the general $\cS_{k}(U)$.

Another motivation, beyond generality, to study $\cS_{k}(U)$ for $U$ other than $\T^n$ comes from the recent preprint of Beck and Everett \cite{beck2026lonely}. They showed that, for each $n$, there exists a finite set $\Lambda$ such that if $\vecb{v}$ is a counterexample to the Lonely Runner Conjecture, then $\langle\vecb{v}\rangle_\R\subseteq H_{\vecb{\lambda}}$ for some $\vecb{\lambda}\in\Lambda$. In other words, if we can show that $\max\cS_{1}(H_{\vecb{\lambda}})\leq\frac12-\frac{1}{n+1}$ for the finitely many $\vecb{\lambda}$ in $\Lambda$, then the Lonely Runner Conjecture holds for $n$. Hence $H_{\vecb{\lambda}}$ is an example of a torus other than $\T^n$ whose relative Lonely Runner spectrum is meaningful to study.

Note that the sequence $\frac12-\frac{s}{ns+1}$ and the sequence in equation \eqref{eq:fansun} are highly arithmetically structured. According to the \emph{structure theorem} of \cite{jain_relative_2026}, this is no coincidence: up to finite symmetric difference every $\cS_1(U)$ for a subtorus with $\dim U=2$ has a similar form. Denote a one-sided arithmetic progression by \[\Prog(\alpha,\beta)=\{\alpha s+\beta:s\in\Z_{\geq0}\}.\]
\begin{theorem}[{\cite[Theorem 1.1]{jain_relative_2026}}]
\label{thm:jainkravitz}
    Let $n\geq2$ and let $U\subseteq\T^n$ be a proper subtorus with $\dim U=2$. There is a finite list $(\alpha_1,\beta_1),\ldots,(\alpha_R,\beta_R)\in\Q_{>0}^2$ such that there is a finite symmetric difference between $\cS_1(U)$ and the set
    \[\bigcup_{r=1}^R\left(D(U)+\frac1{\Prog(\alpha_r,\beta_r)}\right).\]
\end{theorem} 

As the proof of this theorem is constructive, the authors of \cite{jain_relative_2026} were then able to compute some segment of $\cS_1(\T^n)$, up to finite symmetric difference, by enumerating all two-dimensional subtori $U$ with specific values of $D(U)$ and finding all progressions arising from each such $U$. For example, they showed \cite[Theorem 1.3]{jain_relative_2026} that $\cS_1(\T^4)\cap(1/4,1/2]$ and $\frac{1}{4}+\frac{1}{4\Prog(2,3)}$ have finite symmetric difference.

It was then suggested that the analogous statement should remain true for any arbitrary proper subtorus $U\subseteq\T^n$ when the Lonely Runner spectrum is evaluated at codimension one; our main result is this generalization.
\begin{theorem}
\label{thm:main}
    Let $n\geq1$ and let $U\subseteq\T^n$ be a proper subtorus with $\dim U=m$. There is a finite list $(\alpha_1,\beta_1),\ldots,(\alpha_R,\beta_R)\in\Q_{>0}^2$ such that there is a finite symmetric difference between $\cS_{m-1}(U)$ and the set
    \[\bigcup_{r=1}^R\left(D(U)+\frac1{\Prog(\alpha_r,\beta_r)}\right).\]
\end{theorem}

Our proof is similarly constructive, so it would be possible to partially compute new segments of the spectra $\cS_k(\T^n)$ by following the proof of Theorem \ref{thm:main} as in \cite[Sections 3--5]{jain_relative_2026}, but this has not been attempted here. It is unclear whether the computational method arising from our proof improves over Jain and Kravitz's for the case $k=1$, but it can be applied to other values of $k$, so for example it is now possible to compute segments of $\cS_2(\T^4)$ if one wishes. Indeed, it can be applied to an entire class of functions described in Subsection \ref{subsec:rpa}.

As a consequence of Theorem \ref{thm:main}, we can give two descriptions of the classical spectra $\cS_k(\T^n)$. The first concerns the classical spectra in a segment above an accumulation point, and is simply a generalization of \cite[Corollary 1.2]{jain_relative_2026}.
\begin{corollary}
\label{corol:band}
    Let $1\leq k<n$ and let $d\in\acc(\cS_{k}(\T^n))\setminus\cS_{k+2}(\T^n)$ (applying the convention $\cS_{n+1}(\T^n)=\emptyset$). There exist $\varepsilon>0$ and a finite list $(\alpha_1,\beta_1),\ldots,(\alpha_R,\beta_R)\in\Q_{>0}^2$ such that there is a finite symmetric difference between $\cS_{k}(\T^n)\cap(d,d+\varepsilon)$ and the set
    \[\bigcup_{r=1}^R\left(d+\frac1{\Prog(\alpha_r,\beta_r)}\right)\cap(d,d+\varepsilon).\]
\end{corollary}

For completeness, we reproduce the proof from \cite{jain_relative_2026}. From $\acc(\cS_{k+1}(\T^n))\subseteq\cS_{k+2}(\T^n)$, any $d$ above has $\varepsilon>0$ such that $\cS_{k+1}(\T^n)\cap(d-\varepsilon,d+\varepsilon)=\{d\}$. By \cite[Lemmas 4.2 and 4.4]{giri_structure_2026}, there are finitely many proper subtori $U_1,\ldots,U_S$ with $\dim U_s\geq k+1$ such that every $k$-dimensional proper subtorus $T$ (apart from finitely many exceptions) is contained and $\varepsilon$-dense in some $U_s$. Then, for any such $T$ where $D(T)\in(d,d+\varepsilon)$, we have $D(T)-\varepsilon\leq D(U_s)\leq D(T)$ which implies $D(U_s)=d$. If $\dim U_s\geq k+2$, then $d\in\cS_{\dim U_s}(\T^n)\subseteq \cS_{k+2}(\T^n)$ which is a contradiction. Hence $\dim U_s=k+1$ and we apply Theorem \ref{thm:main} on $U_s$. 

The second description concerns the whole  $\cS_k(\T^n)$, and is novel. By choosing the progressions $\Prog(\alpha_r,\beta_r)$ carefully, we can characterize $\cS_k(\T^n)$ as a union of reciprocals of arithmetic progressions, with the only exception being $0$.
\begin{theorem}
\label{thm:arithstruc}
    Let $1\leq k\leq n$. Let $\cU$ be the set of subtori $U\subseteq\T^n$ with $\dim U=k+1$ and $D(U)=0$ (which in particular guarantees properness). There are $(\alpha_{U,1},\beta_{U,1}),\ldots,(\alpha_{U,R_U},\beta_{U,R_U})\in\Q_{>0}^2$ (where each $R_U$ is finite) for each subtorus $U\in\cU$ such that
    \[\cS_k(\T^n)\setminus\{0\}=\bigcup_{U\in\cU}\left(\bigcup_{r=1}^{R_U}\frac1{\Prog(\alpha_{U,r},\beta_{U,r})}\right),\] and moreover, \[\bigcup_{r=1}^{R_U}\frac1{\Prog(\alpha_{U,r},\beta_{U,r})}\subseteq\cS_k(U).\]
\end{theorem}

The proof of this result is more involved and is delayed to Section \ref{sec:classical}, but the following are its implications. First, if a reduced fraction $p/q$ appears in $\cS_k(\T^n)\setminus\{0\}$, then there are infinitely many fractions with numerator $p$ in the spectrum. For instance, the numerator $1$ (and only $1$) appears infinitely often in $\cS_{n-1}(\T^{n})$. 

Moreover, Theorem \ref{thm:arithstruc} suggests an approach to tackle Conjecture \ref{conj:spectra}: if we can design a test that checks whether, for each fixed $(\alpha,\beta)\in\Q_{>0}^2$, we have $1/\Prog(\alpha,\beta)\subseteq\cS_k(\T^n)$, and show that it gives the same result for every $(n,k)$ with the same codimension, we can prove the conjecture. In fact, it should not be difficult to prove that for each $\alpha\in\Q_{>0}$, \[\exists\beta\in\Q_{>0}:\frac{1}{\Prog(\alpha,\beta)}\subseteq\cS_k(\T^n)\iff\exists\beta\in\Q_{>0}:\frac{1}{\Prog(\alpha,\beta)}\subseteq\cS_1(\T^{n-k+1}).\] The difficulty would lie in showing that, for each fixed $\alpha$, the set admissible values of $\beta$ is the same for every $k$.

\subsection{Rational piecewise-affine functions} 
\label{subsec:rpa}
Despite the similarity of the results, our proof of Theorem \ref{thm:main} differs from the proof of Theorem \ref{thm:jainkravitz} in \cite{jain_relative_2026}. The latter exploited the structure of the Lonely Runner problem by observing that on $\T^n$, $D$ is minimized on one of the hyperplanes $\{x_i=\pm x_j\}$ where $i\neq j$. Hence, it suffices to check only the intersections between these hyperplanes and the 2-dimensional $U$, which are one-dimensional subgroups of $\T^n$, which again have rigid arithmetic structure.

On the other hand, our proof does not use the specific properties of the function $D$. Instead, we consider the restriction $D|_U$ under the parent subtorus $U$, which we can interpret as a function from $\T^m$ where $m=\dim U$ to $\R$. Then we study a general class of functions to which this function belongs. Let $\pi:\R^n\rightarrow\T^n$ denote the standard projection onto the torus (which depends implicitly on the dimension of the input).

\begin{definition}
\label{def:rpa}
A function $f:\T^m\to\R$ is \emph{rational piecewise-affine}
if there exist full-dimensional closed polytopes
$P_1,\ldots,P_I\subseteq[0,1]^m$ with vertices in $\Q^m$ such that:
\begin{enumerate}
\item $[0,1]^m=\bigcup_{i=1}^I P_i$,
\item the interiors of the $P_i$ are pairwise disjoint, and
\item for each $i$, there exist $\vecb{a}_i\in\Q^m$ and $b_i\in\Q$
such that
\[
f(\pi(\vecb{y}))=\vecb{a}_i\cdot\vecb{y}+b_i
\qquad
\]
for every $\vecb{y}\in P_i$.
\end{enumerate}
We will occasionally refer to \emph{piecewise-affine} functions, which satisfy the above condition when every instance of $\Q$ is replaced by $\R$.
\end{definition}

Given such $P_1,\ldots,P_I$ and $\vecb{a}_i,b_i$ for each $i$, $f$ is well-defined as a function on $\T^m$ if and only if the values of $f$ on $P_i$ and $P_j$ agree on the intersection $P_i\cap P_j$ and the points on the boundary of $[0,1]^m$ which are identified in the torus quotient. If $f$ is well-defined, this compatibility condition guarantees that $f$ is continuous. 

Rational piecewise-affine functions have been described and studied more abstractly; see, e.g.\ \cite[Definition 2.6.6]{burgos_2014}. However, the question of approximating their minima seems to be completely unexplored. We show that the minima of $f$ on codimension-one subtori follow a highly arithmetic structure. For a continuous $f:\T^m\rightarrow\R$, define \[\cV_f=\{\min_{\vecb{x}\in T} f(\vecb{x}):T\subseteq\T^m\text{ is a }(m-1)\text{-dimensional subtorus}\}.\] Note that $T$ is not required to be proper. 

\begin{theorem}
\label{thm:rpa}
    Let $m\geq1$ be an integer and $f:\T^m\rightarrow\R$ be rational piecewise-affine. There exist $(\alpha_1,\beta_1),\ldots,(\alpha_R,\beta_R)\in\Q_{>0}^2$ (possibly $R=0$) such that there is a finite symmetric difference between $\cV_f$ and the set
    \[\bigcup_{r=1}^R\left(\min_{\vecb{x}\in \T^m}f(\vecb{x})+\frac1{\Prog(\alpha_r,\beta_r)}\right).\] 
\end{theorem}

Theorem \ref{thm:main} is then a direct corollary. Moreover, the more general Theorem \ref{thm:rpa} applies to other interesting functions that can be transformed to a rational piecewise-affine function on $\T^m$. For instance, the \emph{shifted} Lonely Runner Conjecture\footnote{The conjecture of the maximum loneliness being $\geq1/(n+1)$ is false for $5\leq n\leq 16$ \cite{blanco2026coloopless} and $n\geq95$ \cite{poliakova2026more}, and probably false for $16<n<95$ as well. It is still interesting nonetheless to see what the true minimum is.} \cite{alcantara2026covering} asks for the maximum loneliness when each runner is allowed to start at a different point. Fixing the runners' shifts to a rational point $\vecb{c}\in\T^n$
is equivalent to replacing $(1/2,\ldots,1/2)$ in the definition of $D$ by $\vecb{c}-(1/2,\ldots,1/2)$. One can then analogously define the $\vecb{c}$-shifted relative spectra. Theorem \ref{thm:rpa} shows that each codimension-one $\vecb{c}$-shifted relative spectrum has the same structure theorem as the unshifted analogue. 
%\NK{how does this work?  for the definition of the spectrum, are you choosing the ``worst'' shift for each set of speeds?  I agree that your theorem gives something if you fix the shifts in advance and then vary the speeds, but this isn't the same as what sLRC is asking}

\subsection{Preliminaries}
The paper is organized as follows. In Section \ref{sec:codimone}, we prove the codimension-one arithmetic structure Theorem \ref{thm:main}, via the more general Theorem \ref{thm:rpa}. In Section \ref{sec:classical}, we prove the theorem for classical relative spectra, Theorem \ref{thm:arithstruc}.

We fix some notations here. We use boldface italic for vectors and boldface upright for matrices. We use light letters for coordinates of a vector, e.g.\ $\vecb{x}=(x_1,\ldots,x_m)$, and we write $\gcd(\vecb{x})$ for the (nonnegative) gcd of $x_1,\ldots,x_m$. For a function $f:\T^m\rightarrow\R$ and a nonempty and closed subset $H\subseteq\T^m$, we write $d_f(H)=\min_{\vecb{x}\in H}f(\vecb{x})$ and $\delta_f=d_f(\T^m)$.

\section{Codimension-one relative spectra}
\label{sec:codimone}
Assuming Theorem \ref{thm:rpa}, the main Theorem \ref{thm:main} follows once we interpret $D$ on $U$ as a rational piecewise-affine function; we will precisely formulate this in the last subsection. The rest of this section is devoted to the more difficult task of proving Theorem \ref{thm:rpa}.

\subsection{Proof sketch of Theorem \ref{thm:rpa}}
Fix a rational piecewise-affine function $f:\T^m\rightarrow\R$. The codimension-one subtori are $H_{\vecb{\lambda}}$, where $\vecb{\lambda}\in\Z^m$ is primitive. Then our task is to classify the values of $d_f(H_{\vecb{\lambda}})$.

Let $P_1,\ldots,P_I\subseteq[0,1]^m$ be the polytopes in Definition \ref{def:rpa} for $f$. Consider the union of the boundaries\footnote{$B$ is slightly different from the set of non-differentiability $B'=\{\vecb{x}\in\T^m:f\text{ is not differentiable at }\vecb{x}\}$. Perhaps $B'$ is the more intuitive object and indeed in most cases we can make $B=B'$, but we take the boundaries of $P_i$ in $[0,1]^m$ rather than in $\T^m$ to prevent degenerate cases e.g.\,when $f$ is constant and $P_1=[0,1]^m$. In this case, $B$ is the union of the hyperplanes $\{x_i=0\}$, but $B'=\emptyset$ which of course invalidates the rest of our arguments.} $B=\bigcup_{i=1}^I\partial P_i$ which has finitely many edges, which we extend to lines $\ell_1,\ldots,\ell_J$ in $\R^m$. Let $C_j=\pi(\ell_j)$. Due essentially to the fact that an optimum of a linear function on a polytope is attained at a vertex, $d_f(H_{\vecb{\lambda}})$ can be shown to be attained at one of the $C_j$. Note that we can write
$C_j=\vecb{z}_j+\langle \vecb{v}_j\rangle_\R$ with $\vecb{z}_j$ rational and $\vecb{v}_j\in\Z^m$ primitive. It follows that
\begin{equation}
    \label{eq:approxtrop}
    d_f(H_{\vecb{\lambda}})=\min_{j:\vecb{\lambda}\cdot \vecb{v}_j\ne0}d_f(C_j\cap H_{\vecb{\lambda}}).
\end{equation}

For each $j$, since we force that $C_j$ and $H_{\vecb{\lambda}}$ are transverse in equation \eqref{eq:approxtrop}, the intersection $C_j\cap H_{\vecb{\lambda}}$ contains finitely many rational points, equally spaced on $C_j$. We approximate the minima on these points in a similar fashion to \cite[Lemma 2.5]{jain_relative_2026} but with more generality: there are integers $M_j,B_j\geq1$, and constants $\gamma_{j,\vecb{\rho},\sigma}\in\Q_{\geq0}$ indexed by $\vecb{\rho}\in(\Z/M_j\Z)^m,\sigma\in\{\pm\}$, such that, for
\[\vecb{\lambda}\equiv\vecb{\rho}\pmod{M_j},\qquad
 \operatorname{sgn}(\vecb{\lambda}\cdot\vecb{v}_j)=\sigma,\qquad
 |\vecb{\lambda}\cdot\vecb{v}_j|>B_j\]
 we have
\begin{equation}
\label{eq:approxmin}
    d_f(C_j\cap H_{\vecb{\lambda}})
 =d_f(C_j)+\frac{\gamma_{j,\vecb{\rho},\sigma}}{|\vecb{\lambda}\cdot \vecb{v}_j|}.
\end{equation}
Since $f$ is Lipschitz and every point of $\T^m$ is within
$1/(2\|\vecb{\lambda}\|_2)$ of $H_{\vecb{\lambda}}$, there is a uniform bound
\[
 0\leq d_f(H_{\vecb{\lambda}})-\delta_f\leq\frac{A}{\|\vecb{\lambda}\|_2}.
\]
Combining this with equation \eqref{eq:approxmin}, when $\vecb{\lambda}$ is sufficiently large we can ignore all $j$ with $d_f(C_j)>\delta_f$. Relabel the remaining $j$ by $\{1,\ldots,J'\}$.

As a consequence, there is an integer $M>0$ such that we can partition the parameter space for $\vecb{\lambda}$, namely $\Z^m$ with $\vecb{0}$ and non-primitive vectors removed, into finitely many parts according to residues modulo $M$ and some rational linear inequalities (depending on $\vecb{\lambda}\bmod{M}$ and describing the $B_j$ thresholds and comparisons between the various $C_j$'s) so that the following holds. For $\vecb{\lambda}$ sufficiently large, on each part $P$, either $d_f(H_{\vecb{\lambda}})=\delta_f$ or there is $\gamma_P\in\Q_{>0}$ such that
\[
d_f(H_{\vecb{\lambda}})=\delta_f+\frac{\gamma_P}{\ell_P(\vecb{\lambda})}
\]
for some linear form $\ell_P$ with integer coefficients. 

Finally, we employ a number theoretic argument to show that the values
of $\ell_P(\vecb{\lambda})$ on $P$ have finite symmetric difference
with a finite union of arithmetic progressions; the aim is to show that all values in the union can be achieved by \emph{primitive} choices of $\vecb{\lambda}$. As mentioned in  \cite{jain_relative_2026}, this is the part of Jain and Kravitz's proof whose generalization is least obvious. Consequently, for each such $\Prog(a_{P,i},b_{P,i})$ corresponding to $\ell_P(\vecb{\lambda})$, we have
\[
 \delta_f+\frac{\gamma_{P}}{\Prog(a_{P,i},b_{P,i})}
 =\delta_f+\frac1{\Prog(a_{P,i}/\gamma_{P},b_{P,i}/\gamma_{P})}
\]
and the conclusion of Theorem \ref{thm:rpa} follows.

\subsection{Reducing minima to $C_j$} 
Our first claim is equation \eqref{eq:approxtrop}, that the minimum of $f$ on $H_{\vecb{\lambda}}$ is equal to the minimum restricted to the lines $C_j=\vecb{z}_j+\langle \vecb{v}_j\rangle_\R$ intersected with $H_{\vecb{\lambda}}$.

\begin{lemma}
\label{lem:edge_reduction}
Let $m\geq2$ and let $f:\T^m\rightarrow\R$ be piecewise-affine. Let $C_j$ be as previously defined. For every $\vecb{\lambda}\in\Z^m\setminus\{\vecb{0}\}$, we have
\begin{equation}
\label{eq:dfHlamb}
    d_f(H_{\vecb{\lambda}})=\min_{j:\vecb{\lambda}\cdot\vecb{v}_j\ne0}d_f(C_j\cap H_{\vecb{\lambda}}).
\end{equation}
\end{lemma}

\begin{proof}
Choose $\vecb{x}\in H_{\vecb{\lambda}}$ with
$f(\vecb{x})=d_f(H_{\vecb{\lambda}})$, and let
$\vecb{y}\in P_i\subseteq[0,1]^m$ be such that $\pi(\vecb{y})=\vecb{x}$. The polytope
\[
 Q=P_i\cap\{\vecb{u}\in\R^m:\vecb{\lambda}\cdot\vecb{u}=\vecb{\lambda}\cdot\vecb{y}\}
\]
is nonempty because $\vecb{y}\in Q$. Since $f$ is affine on $P_i\supseteq Q$, it attains its minimum on $Q$ at a vertex $\vecb{w}$ of $Q$, so $f(\pi(\vecb{w}))\leq f(\vecb{x})$. However, $f(\vecb{x})$ is minimal, so $f(\pi(\vecb{w}))=f(\vecb{x})$.

Let $F$ be the smallest face of $P_i$ containing $\vecb{w}$.
Then $\vecb{w}$ lies in the relative interior of $F$; otherwise $\vecb{w}$ lies on a smaller face of $F$, contradicting minimality.
If $\dim F\geq2$, since the hyperplane
$\{\vecb{u}\in\R^m:\vecb{\lambda}\cdot\vecb{u}
=\vecb{\lambda}\cdot\vecb{y}\}$ has dimension $m-1$, its intersection with the affine span of $F$ has dimension $\geq1$. Hence $Q$ contains a line through
$\vecb{w}$ parallel to the affine span of $F$.  A sufficiently
short segment of this line, extending on both sides of
$\vecb{w}$, lies in $F\cap Q$. This contradicts that
$\vecb{w}$ is a vertex of $Q$. Hence $\dim F\leq1$.

If $\dim F=1$, then $F$ is an edge of $P_i$ contained in some line $\ell_j$ which is transverse to $H_{\vecb{\lambda}}$ (meaning $\vecb{\lambda}\cdot\vecb{v}_j\neq0$), otherwise $\vecb{w}$ cannot be a vertex of $Q$. If $\dim F=0$, then $\vecb{w}$ is a vertex of $P_i$. Since $P_i$
is full-dimensional, the directions of its edges incident to
$\vecb{w}$ span $\R^m$, hence at least
one such edge is contained in a line $\ell_j$ with
$\vecb{\lambda}\cdot\vecb{v}_j\ne0$.

Thus, in either case, there is $j$ such that $d_f(C_j\cap H_{\vecb{\lambda}})=f(\pi(\vecb{w}))=f(\vecb{x})$ and $\vecb{\lambda}\cdot\vecb{v}_j\ne0$, proving our assertion.
\end{proof}

\subsection{Approximating minima}
The goal of this subsection is to prove that on any line in the torus, if we sample with increasingly fine rational (one-dimensional) grids, the minimum on the grid is eventually of the form in equation \eqref{eq:approxmin}. This result is already implicit in the following lemma of \cite{jain_relative_2026} that dealt with a function on $\T$, but since we will apply it to a function on $C=\vecb{z}+\langle\vecb{v}\rangle_\R$ we will prove the equivalent Lemma \ref{lem:circle_approximation} assuming Lemma \ref{lem:circle_minimum}. 

\begin{lemma}[{\cite[Lemma 2.5]{jain_relative_2026}}]
\label{lem:circle_minimum}
Let $g:\T\to\R$ be rational piecewise-affine and $b\in\Q$ be fixed. There are integers $L,B\geq1$ and constants
$\gamma_r\in\Q_{\geq0}$ indexed by $r\in\Z/L\Z$ such that,
for every integer $q>B$ with $q\equiv r\pmod{L}$,
\[
 \min_{t\in(b+\Z)/q}g(t)
 =\delta_g+\frac{\gamma_r}{q}.
\]
\end{lemma}

In the case of arbitrary lines in the torus, the parameter that replaces the denominator $q$ is the size of the dot product $\vecb{\lambda}\cdot\vecb{v}_j$.

\begin{lemma}
\label{lem:circle_approximation}
Let $f:\T^m\rightarrow\R$ be rational piecewise-affine. Let $C=\vecb{z}+\langle\vecb{v}\rangle_\R$ with $\vecb{z}$ rational and $\vecb{v}\in\Z^m$ primitive. Then there are integers $M,B\geq1$ and constants $\gamma_{\vecb{\rho},\sigma}\in\Q_{\geq0}$ indexed by $\vecb{\rho}\in(\Z/M\Z)^m, \sigma\in\{\pm\}$
such that, whenever $\vecb{\lambda}\in\Z^m$ satisfies
\begin{equation}
\label{eq:manycond}
    \vecb{\lambda}\equiv\vecb{\rho}\pmod M,\qquad
    \operatorname{sgn}(\vecb{\lambda}\cdot\vecb{v})=\sigma,
    \qquad
     |\vecb{\lambda}\cdot\vecb{v}|>B,
\end{equation}
we have
\begin{equation}
\label{eq:approxdenom}
 d_f(C\cap H_{\vecb{\lambda}})
 =d_f(C)+
 \frac{\gamma_{\vecb{\rho},\sigma}}
 {|\vecb{\lambda}\cdot\vecb{v}|}.
\end{equation}
\end{lemma}
\begin{proof}
The map $g:\T\rightarrow\R$ defined by
 $t\mapsto f(\vecb{z}+t\vecb{v})$
is easily seen to be rational piecewise-affine and has $\min g=d_f(C)$. Since $\vecb{z}$ is rational, there is a common denominator $D\in\Z_{>0}$ of every coordinate of $\vecb{z}$. For each $b\in\{0,\frac{1}{D},\ldots,\frac{D-1}{D}\}$, we apply Lemma \ref{lem:circle_minimum} with $g$ and $b$ to obtain integers $L_b,B_b\geq1$ and constants
$c_{b,r}\in\Q_{\geq0}$ indexed by $r\in\Z/L_b\Z$. Then set
\[
\begin{split}
    M&=\operatorname{lcm}(D,L_0,L_{1/D},\ldots,L_{(D-1)/D})\\
    B&=\max(B_0,B_{1/D},\ldots,B_{(D-1)/D}).
\end{split}
\]

Fix a pair $(\vecb{\rho},\sigma)$. We need to pick a corresponding constant $\gamma_{\vecb{\rho},\sigma}$. But before doing that, we pick $b$ as the fractional part $\{-\sigma\vecb{\rho}\cdot\vecb{z}\}$ and $r=(\sigma\vecb{\rho}\cdot\vecb{v}\bmod{L_b})$. At last, pick $\gamma_{\vecb{\rho},\sigma}=c_{b,r}$. 

Now suppose $\vecb{\lambda}$ satisfies \eqref{eq:manycond}. Then put $q=|\vecb{\lambda}\cdot\vecb{v}|=\sigma\vecb{\lambda}\cdot\vecb{v}$. From $\vecb{\lambda}\equiv\vecb{\rho}\pmod M$, we have $\{-\sigma\vecb{\lambda}\cdot\vecb{z}\}=\{-\sigma\vecb{\rho}\cdot\vecb{z}\}=b$, which gives
\begin{align*}
 \pi(\vecb{z}+t\vecb{v})\in H_{\vecb{\lambda}}
 &\iff
 \vecb{\lambda}\cdot\vecb{z}
       +t\vecb{\lambda}\cdot\vecb{v}\in\Z\\
 &\iff qt\in b+\Z.
\end{align*}
Consequently, 
\begin{equation}
    \label{eq:dfchlambda}
    d_f(C\cap H_{\vecb{\lambda}})=\min_{t\in(b+\Z)/q}g(t).
\end{equation}

Finally we want to apply Lemma \ref{lem:circle_minimum} with the chosen $b$ and $q$, but we need to check its conditions. Note that $q>B\geq B_b$, and moreover $q\equiv r\pmod{L_b}$ because  $\vecb{\lambda}\equiv\vecb{\rho}\pmod M$. Hence, we apply the lemma with $b,q$ to obtain
\begin{equation}
    \label{eq:dfCmin}
    \min_{t\in(b+\Z)/q}g(t)=d_f(C)+\frac{c_{b,r}}q=d_f(C)+
 \frac{\gamma_{\vecb{\rho},\sigma}}
 {|\vecb{\lambda}\cdot\vecb{v}|}.
\end{equation}
Combining equations \eqref{eq:dfchlambda} and \eqref{eq:dfCmin} gives the desired equation \eqref{eq:approxdenom}. 
\end{proof}

\subsection{Partitioning the parameter space}
Before we prove the main lemma of this subsection, let us prove a small bound.
\begin{lemma}
\label{lem:lipschitz_bound}
    If $f:\T^m\rightarrow\R$ is $\alpha$-Lipschitz, then
    \begin{equation}
    \label{eq:finite_formulas_bound}
     d_f(H_{\vecb{\lambda}})-\delta_f
     \leq\frac{\alpha}{2\|\vecb{\lambda}\|_2}.
    \end{equation}
\end{lemma}
\begin{proof}
    Let $f$ attain $\delta_f$ at $\vecb{x}$. Let $r$ be the integer closest to $\vecb{\lambda}\cdot\vecb{x}$, so $|r-\vecb{\lambda}\cdot\vecb{x}|\leq1/2$. Take the point
    \[\vecb{y}=\vecb{x}+
     \frac{r-\vecb{\lambda}\cdot\vecb{x}}
     {\|\vecb{\lambda}\|_2^2}\vecb{\lambda},\]
    which has $\vecb{\lambda}\cdot\vecb{y}=r\in\Z$, so $\vecb{y}\in H_{\vecb{\lambda}}$. Lastly, note that $\|\vecb{y}-\vecb{x}\|_2\leq 1/(2\|\vecb{\lambda}\|_2)$.
\end{proof}

The main lemma shows that we can partition $\Z^m$ into finitely many classes by residue class modulo $M$ and rational polyhedra, so that in each set there is a uniform formula for $d_f(H_{\vecb{\lambda}})$ (provided that $\vecb{\lambda}$ is large enough).

\begin{lemma}
\label{lem:finite_formulas}
Let $m\geq2$ and let $f:\T^m\to\R$ be rational piecewise-affine. There are $M\geq1$
and $R>0$ such that for each fixed residue class $\vecb{\rho}\in(\Z/M\Z)^m$ the following holds: There are closed rational polyhedra $Q_1,\ldots,Q_S$ (depending on $\vecb{\rho}$) with the properties that any $\vecb{\lambda}\equiv\vecb{\rho}\pmod{M}$ belongs to (at least) one $Q_s$, and for each $s$ either
\begin{enumerate}
    \item every primitive $\vecb{\lambda}\in Q_s$ with $\vecb{\lambda}\equiv\vecb{\rho}\pmod{M}$ and $\|\vecb{\lambda}\|_2>R$ has $d_f(H_{\vecb{\lambda}})=\delta_f$, or
    \item there is a constant $\gamma_s\in\Q_{>0}$ and a linear form $\ell_s$ positive on $Q_s$ with integer coefficients such that, for every primitive $\vecb{\lambda}\in Q_s$ with $\vecb{\lambda}\equiv\vecb{\rho}\pmod{M}$ and $\|\vecb{\lambda}\|_2>R$, 
    \[
     d_f(H_{\vecb{\lambda}})
     =\delta_f+\frac{\gamma_s}{\ell_s(\vecb{\lambda})}.
    \]
\end{enumerate}
\end{lemma}

\begin{proof}
Apply Lemma \ref{lem:circle_approximation} to every
$C_j=\vecb{z}_j+\langle\vecb{v}_j\rangle_\R$ from
Lemma \ref{lem:edge_reduction}, and let $M_j$ and $B_j$ be the resulting modulus and lower bound respectively (recall that $j\in\{1,\ldots,J\}$). Let $D$ be a common denominator of $\vecb{z}_j$, then choose the common modulus $M=\operatorname{lcm}(D,M_1,\ldots,M_J)$.

Fix $\vecb{\rho}\in(\Z/M\Z)^m$, and write
$a_j(\vecb{x})=\vecb{x}\cdot\vecb{v}_j$. For fixed $j$, each $\vecb{\lambda}\equiv\vecb{\rho}\pmod{M}$ belongs to one of the following $2B_j+3$ possibilities:
\begin{enumerate}
    \item $k-1/2\leq a_j(\vecb{\lambda})\leq k+1/2$, i.e.\,$a_j(\vecb{\lambda})=k$, for some integer $-B_j\leq k\leq B_j$
    \item $a_j(\vecb{\lambda})\geq B_j+1/2$
    \item $a_j(\vecb{\lambda})\leq-B_j-1/2$.
\end{enumerate}
Let these $\sum_{j=1}^J(2B_j+3)$ inequalities, called \emph{Class I} inequalities, define finitely many closed
rational polyhedra $T_1,\ldots,T_W$.

Now fix $T_w$ and let $\vecb{\lambda}\in T_w$. Consider the value of $d_f(C_j\cap H_{\vecb{\lambda}})$ in the right hand side of equation \eqref{eq:dfHlamb} (in Lemma \ref{lem:edge_reduction}). If $a_j=0$, this $C_j$ is ignored when we take the minimum. If $a_j=k$ is fixed where $|k|\in\{1,\ldots,B_j\}$,
then
\[
 \pi(\vecb{z}_j+t\vecb{v}_j)\in H_{\vecb{\lambda}}
 \iff kt+\vecb{\rho}\cdot\vecb{z}_j\in\Z.
\]
This set of points depends only on $k$, so $d_f(C_j\cap H_{\vecb{\lambda}})$ is a rational number depending only on $k$, say $\delta_{j,k}$.
Otherwise, $|a_j|>B_j$, so Lemma \ref{lem:circle_approximation} gives $\gamma_j\in\Q_{\geq0}$, and
$q_j=\sigma_j a_j$ where $\sigma_j\in\{\pm\}$ is chosen so that $q_j$ is positive on $T_w$, such that
\begin{equation}
    \label{eq:dflargeest}
    d_f(C_j\cap H_{\vecb{\lambda}})
 =d_f(C_j)+\frac{\gamma_j}{q_j(\vecb{\lambda})}.
\end{equation}

Since $f$ is clearly Lipschitz, by Lemma \ref{lem:lipschitz_bound}, we can choose $R_{\vecb{\rho}}$ large enough so that if $\|\vecb{\lambda}\|_2>R_{\vecb{\rho}}$ then any $C_j$ with $d_f(C_j)>\delta_f$ or $\delta_{j,k}>\delta_f$ cannot have minimal $d_f(C_j\cap H_{\vecb{\lambda}})$. We can make this $R$ uniform by taking the maximum of all $R_{\vecb{\rho}}$. Henceforth also assume $\|\vecb{\lambda}\|_2>R$.

If there is $j$ such that $\delta_{j,k}=\delta_f$, then in $T_w$, we have $d_f(H_{\vecb{\lambda}})=\delta_f$; this is conclusion (1) in the lemma. 

If there is no such $j$, there are finitely many $j$ such that $d_f(C_j)=\delta_f$ (and at least one, because otherwise Lemma \ref{lem:edge_reduction} fails), which we name $j_1,\ldots,j_p$. From equation \eqref{eq:dflargeest}, to find the smallest $d_f(C_j\cap H_{\vecb{\lambda}})$ we need to find the smallest $\gamma_{j_1}/q_{j_1}(\vecb{\lambda}),\ldots,\gamma_{j_p}/q_{j_p}(\vecb{\lambda})$. Observe that
\[
 \frac{\gamma_i}{q_i(\vecb{\lambda})}
 \leq\frac{\gamma_j}{q_j(\vecb{\lambda})}
 \iff
 \gamma_iq_j(\vecb{\lambda})
 \leq\gamma_jq_i(\vecb{\lambda}),
\]
for $i,j\in\{j_1,\ldots,j_p\}$. These are rational linear inequalities, which we call \emph{Class II} inequalities. Hence we subdivide each such $T_w$ by these Class II inequalities, and we have conclusion (2) for each subdivided polyhedron.

The closed rational polyhedra $Q_1,\ldots,Q_S$ in the lemma statement are the polyhedra defined by the Class I and Class II inequalities.
\end{proof}

\subsection{Obtaining arithmetic progressions}
For vectors $\vecb{b},\vecb{u}\in\Z^m$, we denote the gcd of the minors of $(\vecb{b}\mid\vecb{u})$ by
\begin{equation}
    \label{eq:gcdminor}
    \Delta(\vecb{b},\vecb{u})=\gcd_{i<j}(b_iu_j-b_ju_i).
\end{equation} 
We make an elementary observation that $\Delta(\vecb{b},\vecb{u})$ controls the gcd of a vector in the arithmetic progression with common difference $\vecb{u}$ and shift $\vecb{b}$.
\begin{lemma}
\label{lem:gcd}
    Let $\vecb{b},\vecb{u}\in\Z^m$ be linearly independent, so $\Delta(\vecb{b},\vecb{u})\neq0$. For every $a\in\Z$, $\gcd(\vecb{b}+a\vecb{u})$ divides $\Delta(\vecb{b},\vecb{u})$.
\end{lemma}

\begin{proof}
    As $\gcd(\vecb{b}+a\vecb{u})$ divides every coordinate of $\vecb{b}+a\vecb{u}$, it divides
    \[(b_i+au_i)u_j-(b_j+au_j)u_i=b_i u_j-b_j u_i\]
    for every $i<j$. It therefore divides their gcd $\Delta(\vecb{b},\vecb{u})$.
\end{proof}

Then, we can show that the linear function $\ell$ of each class from the previous subsection contributes a finite union of arithmetic progressions, up to finite symmetric difference. 

\begin{lemma}
\label{lem:arithmetic}
Let $P\subseteq\R^m$ be a rational polyhedron, and let $\ell$ be a linear form positive on $P$ with integer coefficients. Fix integer $M\geq1$ and $\vecb{\rho}\in(\Z/M\Z)^m$. Then, the set
\[
 E=\{\ell(\vecb{\lambda}):
 \vecb{\lambda}\in P,\ \vecb{\lambda}\equiv\vecb{\rho},
 \ \gcd(\vecb{\lambda})=1\}
\]
has finite symmetric difference with a finite (possibly empty) union of arithmetic progressions $\Prog(a,b)$ with $a,b\in\Z_{>0}$.
\end{lemma}
The idea of the proof is to observe that, if we do not restrict $\gcd(\vecb{\lambda})=1$ for $E$, then the conclusion follows immediately by choosing a ray for each progression. Hence, it suffices to show that, for each $a,b$ in our finite union above, every sufficiently large term of $\Prog(a,b)$ is produced by a primitive $\vecb{\lambda}$. Indexing the terms by $n$, this leads us to consider the level sets of $\ell$ given by \[\{\vecb{c}_n+Mt\vecb{v}:t\in\Z\cap[-n,n]\}\] (the constants will be defined in the proof). The number theory part is then to show that for $n$ sufficiently large, at least one of the vectors $\vecb{c}_n+Mt\vecb{v}$ is primitive. However, we can only produce $\vecb{c}_n$ and $\vecb{v}$ for these level sets when the dimension of the unbounded part of $P$ is at least $2$, therefore we give a separate argument for when this dimension is $1$.

\begin{proof}[Proof of Lemma \ref{lem:arithmetic}]
If $\ell(P)$ is bounded, $E$ is finite so the statement is trivially true. If $\gcd(M,\rho_1,\ldots,\rho_m)>1$, $E=\emptyset$. So assume that
$\ell(P)$ is unbounded and
$\gcd(M,\rho_1,\ldots,\rho_m)=1$.

We can decompose (see e.g.\ \cite[Corollary 7.1b]{schrijver1998theory})
\begin{equation}
    \label{eq:decompose}
    P=K+\langle\vecb{u}_1,\ldots,\vecb{u}_I\rangle_{\R_{\geq0}},
\end{equation}
where $K$ is a bounded polytope and the $\vecb{u}_i$ are primitive
integer vectors. By linearity note that if $\vecb{b}\in K$, then $0<\ell(\vecb{b}+t\vecb{u_1})=\ell(\vecb{b})+t\ell(\vecb{u_1})$ for every $t\in\R_{>0}$, hence $\ell(\vecb{u}_i)\geq0$. Since $\ell(P)$ is unbounded, we may assume without loss that $\ell(\vecb{u}_1)>0$.

The first case is when all $\vecb{u}_i$ are parallel. Then $\vecb{u}_i$ is not a negative scalar multiple of $\vecb{u}_j$, as $\ell(P)$ is unbounded and positive, hence $P=K+\langle\vecb{u}_1\rangle_{\R_{\geq0}}$. 
%\NK{can't you always assume that $I=m$ and that the $u$'s are linearly independent?  then I suppose that this paragraph would correspond to the case $m=1$}
Let $B=\left(K+[0,M]\vecb{u}_1\right)\cap(\vecb{\rho}+M\mathbb Z^m)$, which is bounded hence finite. Then \[P\cap(\vecb{\rho}+M\Z^m)=\bigcup_{\vecb{b}\in B}\left(\vecb{b}+M\langle\vecb{u}_1\rangle_{\Z_{\geq0}}\right).\] Fix $\vecb{b}\in B$ and set $\Delta=\Delta(\vecb{b},\vecb{u}_1)$.
If $\Delta=0$, then $\vecb{b}$ and $\vecb{u}_1$ are parallel, hence there are at most $2$ primitive vectors in the ray $\vecb{b}+M\langle\vecb{u}_1\rangle_{\Z_{\geq0}}$. If $\Delta\neq0$, by Lemma \ref{lem:gcd},
\[
 \gcd(\vecb{b}+Mt\vecb{u}_1)=1
 \iff
 \gcd(\Delta,b_1+Mt u_{1,1},\ldots,b_m+Mt u_{1,m})=1.
\]
The primitivity of $\vecb{b}+Mt\vecb{u}_1$ depends only on $t\bmod\Delta$. Hence, for fixed $d\in\{0,1,\ldots,\Delta-1\}$, either all of $\vecb{b}+M(\Delta t+d)\vecb{u}_1$ are primitive, or none are. Thus each $\vecb{b}\in B$ contributes either finitely many arithmetic progressions or finitely many values to $E$. Since $B$ is finite, we have the conclusion for $E$.

The second case is when not all $\vecb{u}_i,\vecb{u}_j$ are parallel, so assume that $\vecb{u}_2$ is not parallel to $\vecb{u}_1$. Let
\begin{gather*}
 \vecb{w}=\ell(\vecb{u}_1)\vecb{u}_2-\ell(\vecb{u}_2)\vecb{u}_1
\end{gather*}
and $\vecb{v}=\vecb{w}/\gcd(\vecb{w})$.
Thus $\vecb{v}$ is primitive and $\ell(\vecb{v})=0$.
Choose an integer $N$ large enough that, with
$\vecb{u}=N(\vecb{u}_1+\vecb{u}_2)$,
\begin{equation}
    \label{eq:upmvinclude}
    \vecb{u}+\vecb{v},\ \vecb{u}-\vecb{v}\in\langle\vecb{u}_1,\vecb{u}_2\rangle_{\R_{\geq0}}.
\end{equation}

Now take $A=M\ell(\vecb{u})>0$ and fix $a\in\Z/A\Z$. Suppose that there exists a primitive $\vecb{b}\in P\cap(\vecb{\rho}+M\Z^m)$ such that $\ell(\vecb{b})\equiv a\pmod{A}$. We claim that $\ell$ produces $\Prog(A,a)$ with at most finitely many exceptions. For every integer $n\geq1$, set
\[
 \vecb{c}_n=\vecb{b}+Mn\vecb{u}.
\]
Observe that \[ \vecb{c}_n+Mt\vecb{v}=\vecb{b}+M(n\vecb{u}+t\vecb{v})
 =\vecb{b}+M\left(\frac{n+t}{2}(\vecb{u}+\vecb{v})
  +\frac{n-t}{2}(\vecb{u}-\vecb{v})\right),\]
so $\vecb{c}_n+Mt\vecb{v}\in P$ if the last term belongs in $\langle\vecb{u}_1,\vecb{u}_2\rangle_{\R_{\geq0}}$. By \eqref{eq:upmvinclude}, this happens when $|t|<n$. Moreover, if $t\in\Z$ as well, then \[\vecb{c}_n+Mt\vecb{v}\in P\cap(\vecb{\rho}+M\Z^m).\] Lastly note that $\ell(\vecb{c}_n+Mt\vecb{v})=\ell(\vecb{b})+nA$. Therefore, it suffices to show that when $n$ is large, there is a choice of $t\in\Z,|t|<n$ such that $\vecb{c}_n+Mt\vecb{v}$ is primitive.

Define $\Delta_n=\Delta(\vecb{c}_n,\vecb{v})$.
Since $\ell(\vecb{c}_n)>0$ and $\ell(\vecb{v})=0$,
the vectors $\vecb{c}_n$ and $\vecb{v}$ are linearly independent, so $\Delta_n>0$. Each minor grows linearly in $n$, hence there is a constant $C>0$ such that $1\leq \Delta_n\leq Cn$.

Again by Lemma \ref{lem:gcd}, any prime 
$p\mid\gcd(\vecb{c}_n+Mt\vecb{v})$ divides $\Delta_n$. Clearly, such $p\nmid M$, otherwise $\vecb{b}$ is not primitive. For each prime $p\mid\Delta_n$ with $p\nmid M$, the relations $c_{n,i}v_j-c_{n,j}v_i\equiv0\pmod{p}$ for all $i<j$ implies that modulo $p$, $\vecb{c}_n$ and $\vecb{v}$ are scalar multiples. Since $\vecb{v}\not\equiv\vecb{0}\pmod{p}$, there is only one choice of $t\bmod{p}$ that makes $\vecb{c}_n+Mt\vecb{v}\equiv0\pmod{p}$, so only one choice of $t\bmod{p}$ is excluded.

Let $\mathcal{P}_n$ be the set of primes dividing $\Delta_n$ but not $M$ and let $\omega(\Delta_n)$ be the number of primes dividing $\Delta_n$, so $\#\mathcal{P}_n\leq\omega(\Delta_n)$. Applying inclusion--exclusion in the standard way with the sets $\{t\in\Z\cap[-n,n]: p_1\cdots p_k\mid\gcd(\vecb{c}_n+Mt\vecb{v})\}$ where $p_1,\ldots,p_k\in\mathcal{P}_n$, each of which has size $(p_1\cdots p_k)^{-1}(2n+1)$ within error $<1$, we get
\begin{equation}
\label{eq:inclexcl}
    \begin{split}
        \#\{t\in\Z\cap[-n,n]:\gcd(\vecb{c}_n+Mt\vecb{v})=1\}&>
 (2n+1)\prod_{\substack{p\in\mathcal{P}_n}}
 \left(1-\frac1p\right)-2^{\#\mathcal{P}_n}\\
 &\geq (2n+1)2^{-\#\mathcal{P}_n}-2^{\#\mathcal{P}_n}\\
 &\geq (2n+1)2^{-\omega(\Delta_n)}-2^{\omega(\Delta_n)}
    \end{split}
\end{equation}
The last term is positive if $4^{\omega(\Delta_n)}<2n+1$. This holds for sufficiently large $n$, because there are 6 primes below $17$, so
\[
 4^{\omega(\Delta_n)}\leq4^6\sqrt{\Delta_n}\leq4^6\sqrt{Cn}=o(n).
\]
Since there are finitely many residue classes modulo $A$, this proves our claim. %\NK{this is a slick argument.  note that it's non-constructive on bounding the number/size of the finitely many exceptions, since it doesn't give a bound on the size of $\vecb{b}$.  did you also try adapting the more overtly sieve-theoretic argument from my paper with Vanshika?  it would be interesting to remark on whether or not that can also be made to work}
\end{proof}

We make two remarks on this proof:
\begin{enumerate}
    \item Strictly speaking, the proof is non-constructive (on bounding the number of the finitely many exceptions) because, for each $a$, we choose any $\vecb{b}\in P\cap(\vecb{\rho}+M\Z^m)$ such that $\ell(\vecb{b})\equiv a\pmod{A}$ provided that one exists. However, we can bound the size of $\vecb{b}$ that we have to check because if $\vecb{a}=\vecb{b}\bmod{AM}$, then $\vecb{a}\in P\cap(\vecb{\rho}+M\Z^m)$ and $\ell(\vecb{a})\equiv a\pmod{A}$. By the decomposition \eqref{eq:decompose}, we then only have to check $\vecb{b}$ in $P'=K+\langle\vecb{u}_1,\ldots,\vecb{u}_I\rangle_{S}$ where the spanning set $S$ is the real interval $[0,AM]$ (the set $P'$ is the analogue of $B$ in the second case). Hence the proof can be made constructive.
    \item %For each $n$, one can view the line segment $\vecb{c}_n+Mt\vecb{v}$ where $t\in[-n,n]$ as a level set of $\ell$. 
    The estimate \eqref{eq:inclexcl} gives \[\lim_{n\rightarrow\infty}\#\{t\in\Z\cap[-n,n]:\gcd(\vecb{c}_n+Mt\vecb{v})=1\}=\infty.\] In other words, a level set of $\ell$ can be made arbitrarily large. This statement only requires $\gcd(\vecb{c_n},M\vecb{v})=1$, and is essentially the $m$-dimensional analogue of \cite[Lemma 2.7]{jain_relative_2026}. Of course, it could have been proven with more analytic number theory like in \cite{jain_relative_2026} if we want a better numerical bound, but we prefer this elementary approach.
\end{enumerate}

\subsection{Putting everything together}
Finally, we combine Lemmas \ref{lem:finite_formulas} and \ref{lem:arithmetic} to obtain our desired theorems that $\cV_f-\delta_f$ is, up to finite symmetric difference, the finite union of the reciprocals of some arithmetic progressions.

\begin{proof}[Proof of Theorem \ref{thm:rpa}]
If $m=1$, the only $0$-dimensional subtorus is $\{\vecb{0}\}$, so $\cV_f$ has finite symmetric difference with the empty progression. For $m\geq2$, Lemma \ref{lem:finite_formulas} produces an integer $M>0$ and, for each $\vecb{\rho}\in(\Z/M\Z)^m$, polyhedra $Q_{\vecb{\rho},1},\ldots,Q_{\vecb{\rho},S_{\vecb{\rho}}}$ with $d_f(H_{\vecb{\lambda}})=\delta_f$ or $d_f(H_{\vecb{\lambda}})=\delta_f+\gamma_{\vecb{\rho},s}/\ell_{\vecb{\rho},s}(\vecb{\lambda})$ uniformly on the polyhedron provided that $\|\vecb{\lambda}\|_2>R$. Let $\cV_f'$ be the set of values of $d_f(H_{\vecb{\lambda}})$ as described above, but with the restriction $\|\vecb{\lambda}\|_2>R$ removed. As the ball $\|\vecb{\lambda}\|_2\leq R$ is bounded, $\cV_f$ and $\cV_f'$ are equal modulo finite symmetric difference. 

If $Q_{\vecb{\rho},s}$ has $d_f(H_{\vecb{\lambda}})=\delta_f$, it contributes a single value to $\cV_f'$. If $Q_{\vecb{\rho},s}$ has $d_f(H_{\vecb{\lambda}})=\delta_f+\gamma_{\vecb{\rho},s}/\ell_{\vecb{\rho},s}(\vecb{\lambda})$, by Lemma \ref{lem:arithmetic} the values attained by $\ell_{\vecb{\rho},s}(\vecb{\lambda})$ have finite symmetric difference with a finite union of $\Prog(a_{\vecb{\rho},s,i},b_{\vecb{\rho},s,i})$ with $a_{\vecb{\rho},s,i},b_{\vecb{\rho},s,i}\in\Z_{>0}$. Since
\[
 \delta_f+\frac{\gamma_{\vecb{\rho},s}}{\Prog(a_{\vecb{\rho},s,i},b_{\vecb{\rho},s,i})}
 =\delta_f+\frac{1}{\Prog(\alpha_{\vecb{\rho},s,i},\beta_{\vecb{\rho},s,i})},
\]
where $\alpha_{\vecb{\rho},s,i}=a_{\vecb{\rho},s,i}/\gamma_{\vecb{\rho},s}$ and $\beta_{\vecb{\rho},s,i}=b_{\vecb{\rho},s,i}/\gamma_{\vecb{\rho},s}$, $\cV_f'$ has finite symmetric difference with 
\[\bigcup_{\vecb{\rho},s,i}\left(\delta_f+\frac1{\Prog(\alpha_{\vecb{\rho},s,i},\beta_{\vecb{\rho},s,i})}\right)\]
(which is a finite union), and hence so does $\cV_f$.
\end{proof}

To prove Theorem \ref{thm:main} from Theorem \ref{thm:rpa}, we only need to show that under a suitable transformation, the function $D(\vecb{x})=\|\vecb{x}-(1/2,\ldots,1/2)\|_\infty$ is rational piecewise-affine.

\begin{proof}[Proof of Theorem \ref{thm:main}]
It suffices to exhibit an isomorphism $\iota:\T^m\rightarrow U$, since $f=D\circ\iota$ is then rational piecewise-affine, and we can apply Theorem \ref{thm:rpa} to $f$. (Although $\cV_f$ is defined for all $(m-1)$-dimensional subtori which are not necessarily proper, non-proper subtori $T$ only give $\min_{\vecb{x}\in T}f(\vecb{x})=1/2$, so the conclusion still holds.) 

To construct $\iota$, let $V\subseteq \mathbb R^n$ be the $m$-dimensional rational subspace such that $U=\pi(V)$. Then the lattice $V\cap\mathbb Z^n$ has rank $m$, so there is a $\Z$-basis $\vecb{u}_1,\ldots,\vecb{u}_m$. Define $\iota:\T^m\rightarrow\T^n$ to be induced by the matrix $(\vecb{u}_1|\cdots|\vecb{u}_m)$. Then $\iota$ is a homomorphism with $\operatorname{Im}(\iota)=U$. Since the lattice is saturated, we can extend $\vecb{u}_1,\ldots,\vecb{u}_m$ to a $\Z$-basis of $\Z^n$ denoted by  $\vecb{u}_{1},\ldots,\vecb{u}_n$. Set $\mat{P}=(\vecb{u}_1|\cdots|\vecb{u}_n)$. The homomorphism $\kappa:\T^n\rightarrow\T^m$, induced by the first $m$ rows of $\mat{P}^{-1}$, gives $\kappa\circ\iota=I_m$. Therefore, $\iota$ has an inverse $\kappa|_U$, and is thus an isomorphism.
\end{proof}

\section{Classical spectra}
\label{sec:classical}
The goal of this section is to prove Theorem \ref{thm:arithstruc} concerning the classical spectra. Let $1\leq k<n$ be fixed and recall that we defined \[\cU=\{U\subseteq\T^n:\dim U=k+1,D(U)=0\}.\] Consider a fixed subtorus $U_0\in\cU$. By removing finitely many initial terms from each arithmetic progression in Theorem \ref{thm:main}, we can make the union of reciprocals of arithmetic progressions be a subset of $\cS_k(U_0)$. The next lemma says that, for any $d\in\cS_k(U_0)\setminus\{0\}$, there is a different $U\in\cU$ whose spectrum contains the reciprocals of an arithmetic progression, one of which is $d$. We apply this lemma to the values $d$ removed from $\cS_k(U_0)$.

\begin{lemma}
\label{lem:amend}
    Let $1\leq k<n$, and let $T\subseteq\T^n$ be a proper $k$-dimensional subtorus with $d=D(T)>0$. There exist a
    $(k+1)$-dimensional subtorus $U\supseteq T$ with $D(U)=0$ (meaning $U\in\cU$)
    and $\alpha\in\Q_{>0}$ such that
    \[
     \frac1{\Prog(\alpha,1/d)}\subseteq\cS_k(U).
    \]
\end{lemma}
\begin{proof}
Put $\vecb{c}=(1/2,\ldots,1/2)\in\R^n$ and construct $U$ as follows. Let $\vecb{x}\in[0,1]^n$ be a point such that $\pi(\vecb{x})\in T$ and $\|\vecb{x}-\vecb{c}\|_\infty=d$, hence $\vecb{x}-\vecb{c}$ is the distance vector from the centre to $\vecb{x}$. Let $V\subseteq \mathbb{R}^n$ be the subspace such that $T=\pi(V)$. Then let $W$ be the span of $V$ and $\vecb{x}-\vecb{c}$ and $U=\pi(W)$. As $\pi(\vecb{c})\notin T$, $U$ is a $(k+1)$-dimensional subtorus which is proper because $D(U)=0<1/2$.

Now put $m=k+1$. Let $\iota:\T^m\to U$ be an isomorphism induced by the matrix $\mat{A}\in\Z^{m\times n}$ (as constructed in the proof of Theorem \ref{thm:main}), and let $f=D\circ\iota$. Pick
$\vecb{z}\in\frac12\Z^m$ with $\iota(\vecb{z})=\pi(\vecb{c})$. Then there is a primitive $\vecb{\lambda}_0\in\Z^m$ such that $T=\iota(H_{\vecb{\lambda}_0})$. Moreover, $\vecb{\lambda}_0\cdot\vecb{z}\equiv1/2\pmod1$ because $\pi(\vecb{c})\notin T$. Now define \[\mu(\vecb{\lambda})=\min_{\vecb{\lambda}\cdot\vecb{u}\in\frac12+\Z}\|\mat A\vecb{u}\|_\infty.\] The main idea of the proof is to express $d_f(H_{\vecb{\lambda}})$ as $\mu(\vecb{\lambda})$, which has an easy linear form.

First, we claim that $d_f(H_{\vecb{\lambda}})=\mu(\vecb{\lambda})$ when $\vecb{\lambda}\equiv\vecb{\lambda}_0 \pmod{2}$ is sufficiently large. Since $\vecb{z}\in\frac12\Z^m$, we have $\vecb{\lambda}\cdot\vecb{z}\equiv\vecb{\lambda}_0\cdot\vecb{z}\equiv1/2\pmod1$
and therefore $\pi(\vecb{z}+\vecb{u})\in H_{\vecb{\lambda}}$ if and only if $\vecb{\lambda}\cdot\vecb{u}\in\tfrac12+\Z$. As $f(\pi(\vecb{z}+\vecb{u}))
=\min_{\vecb{q}\in\Z^n}
\|\mat A\vecb{u}-\vecb{q}\|_\infty
\leq\|\mat A\vecb{u}\|_\infty
$, we get
\begin{equation}
\label{eq:normupper}
    d_f(H_{\vecb{\lambda}})=\min_{\vecb{\lambda}\cdot\vecb{u}\in\frac12+\Z}f(\pi(\vecb{z}+\vecb{u}))\leq\min_{\vecb{\lambda}\cdot\vecb{u}\in\frac12+\Z}\|\mat A\vecb{u}\|_\infty=\mu(\vecb{\lambda}).
\end{equation}

For the reverse inequality, choose $\delta>0$ sufficiently small that
$f(\pi(\vecb{z}+\vecb{u}))=\|\mat{A}\vecb{u}\|_\infty$ whenever
$\|\vecb{u}\|_2<\delta$.
Since $f$ is continuous and vanishes only at $\pi(\vecb{z})$,
compactness gives $\varepsilon>0$ such that every point $\vecb{x}$ with
$f(\vecb{x})<\varepsilon$ can be written as $\vecb{x}=\pi(\vecb{z}+\vecb{u})$
with $\|\vecb{u}\|_2<\delta$.
By Lemma \ref{lem:lipschitz_bound}, there exists $R>0$ such that
$d_f(H_{\vecb{\lambda}})<\varepsilon$ whenever
$\|\vecb{\lambda}\|_2>R$.
For such $\vecb{\lambda}$, let $\vecb{x}\in H_{\vecb{\lambda}}$ be a
minimizer of $f$, and choose a lift $\vecb{z}+\vecb{u}$ of $x$
with $\|\vecb{u}\|_2<\delta$. Then
\[
 d_f(H_{\vecb{\lambda}})
 =f(\pi(\vecb{z}+\vecb{u}))
 =\|\mat{A}\vecb{u}\|_\infty
 \geq\mu(\vecb{\lambda}),
\]
as desired.

Second, we claim that $d_f(H_{\vecb{\lambda}_0})=\mu(\vecb{\lambda}_0)$ as well. Equation \eqref{eq:normupper} holds for $\vecb{\lambda}=\vecb{\lambda}_0$, so $d_f(H_{\vecb{\lambda}_0})\leq\mu(\vecb{\lambda}_0)$. On the other hand, our choice of $W$ ensures that $\vecb{x}-\vecb{c}\in W$, hence there is $\vecb{u}_0\in\R^m$ such that $\mat A\vecb{u}_0=\vecb{x}-\vecb{c}$.\footnote{We want $\mat A\vecb{u}_0=\vecb{x}-\vecb{c}$ in $\R^n$ not just $\T^n$, so that we can obtain the first equality in \eqref{eq:firstterm}. If this only holds in $\T^n$, then $\|\mat A\vecb{u}_0\|_\infty
\geq\|\vecb{x}-\vecb{c}\|_\infty$ which is insufficient. This is why we need to consider a lift for $\vecb{x}$ that minimizes the real distance to $\vecb{c}$ (but there is nothing special about this lift of $\vecb{c}$; we may as well choose a different lift $(3/2,\ldots,3/2)$ with its corresponding $\vecb{x}$).}
Since $\mat{A}(\vecb{z}+\vecb{u}_0)=\vecb{x}$, we have $\pi(\vecb{z}+\vecb{u}_0)\in H_{\vecb{\lambda}_0}$ and hence
$\vecb{\lambda}_0\cdot\vecb{u}_0\in\frac12+\Z$. Therefore,
\begin{equation}
    \label{eq:firstterm}
    \mu(\vecb{\lambda}_0)
 \leq\|\mat A\vecb{u}_0\|_\infty
 =\|\vecb{x}-\vecb{c}\|_\infty
 =d
 =d_f(H_{\vecb{\lambda}_0}),
\end{equation}
which proves the equality.

Similarly to the first case in the proof of Lemma \ref{lem:arithmetic}, we consider a ray. Set the sequence \[\vecb{\lambda}_s=\vecb{\lambda}_0+sK\vecb{v}\]
for $s\in\Z_{\geq0}$ and $T_s=\iota(H_{\vecb{\lambda}_s})\subseteq U$, with $K$ and $\vecb{v}$ given as follows. Let $P=\{\vecb{x}\in\R^m:\|\mat{A}\vecb{x}\|_\infty\leq1\}$ be a rational polytope, which is full-dimensional and centrally symmetric. Let $\vecb{w}\in P$ be a vertex such that $\vecb{\lambda}_0\cdot\vecb{w}$ is maximized. The rational polyhedron
\[
 C=\{\vecb{\lambda}\in\R^m:
 \vecb{\lambda}\cdot\vecb{w}
 \geq\vecb{\lambda}\cdot\vecb{u}
 \text{ for every }\vecb{u}\in P\}
\] is full-dimensional and in particular $\vecb{\lambda}_0\in C$. Then choose $\vecb{v}\in C\cap\Z^m$ linearly independent of $\vecb{\lambda}_0$. Furthermore, let $\Delta=\Delta(\vecb{\lambda}_0,\vecb{v})$ be the gcd of the minors of $(\vecb{\lambda}_0\mid\vecb{v})$, so $\Delta\neq0$, and pick $K$ a multiple of $2\Delta$ with
$K\|\vecb{v}\|_2>R+\|\vecb{\lambda}_0\|_2$ (so that $\|\vecb{\lambda}_0+sK\vecb{v}\|_2\geq sK\|\vecb{v}\|_2-\|\vecb{\lambda}_0\|_2>R$ for $s\geq1$).

We claim that $\vecb{\lambda}_s$ is primitive. By Lemma \ref{lem:gcd}, $\gcd(\vecb{\lambda}_s)\mid\Delta$, and since $\Delta\mid K$, if $\vecb{\lambda}_s$ is not primitive then a prime $p\mid\gcd(\vecb{\lambda}_s)$ gives $p\mid\gcd(\vecb{\lambda}_0)$, contradicting that $\vecb{\lambda}_0$ is primitive.

It remains to show that $\mu$ indeed has a reciprocal-linear form in $s$. One can check that $\vecb{\lambda}_s\in C$ and
$\vecb{\lambda}_s\equiv\vecb{\lambda}_0\pmod2$. By symmetry and homogeneity, for every $\vecb{\lambda}\in C\setminus\{\vecb{0}\}$, we have
\[
 \mu(\vecb{\lambda})
 =\frac1{2\vecb{\lambda}\cdot\vecb{w}}.
\]
Hence, finally,
\[
 D(T_s)=d_f(H_{\vecb{\lambda}_s})=\mu(H_{\vecb{\lambda}_s})
 =\frac1{2\vecb{\lambda}_s\cdot\vecb{w}}
 =\frac1{d^{-1}+2K(\vecb{v}\cdot\vecb{w})s}
\]
for every $s\in\Z_{\geq0}$.
Taking $\alpha=2K(\vecb{v}\cdot\vecb{w})\in\Q_{>0}$ proves the assertion.
\end{proof}

Theorem \ref{thm:arithstruc} follows easily.

\begin{proof}[Proof of Theorem \ref{thm:arithstruc}]
    If $k=n$, then this is obvious since $\cS_n(\T^n)=\{0\}$ and $\cU=\emptyset$. Now let $k<n$.
    
    Since every $k$-dimensional $T\subseteq\T^n$ is included in some $U\in\cU$ (e.g.\ by the construction in the proof of Lemma \ref{lem:amend}), we have $\cS_k(\T^n)=\bigcup_{U\in\cU}\cS_k(U)$. For each $U\in\cU$, shorten the arithmetic progressions given by Theorem \ref{thm:main} so that we can write $\mathcal{P}_U\subseteq\cS_k(U)$ where \[\mathcal{P}_U=\bigcup_{r=1}^{R_U}\left(D(U)+\frac1{\Prog(\alpha_{U,r},\beta_{U,r})}\right)=\bigcup_{r=1}^{R_U}\left(\frac1{\Prog(\alpha_{U,r},\beta_{U,r})}\right)\] and the set of exceptions $\mathcal{E}_U=\cS_k(U)\setminus\mathcal{P}_U$ is finite. Then let \[\mathcal{E}=\cS_k(\T^n)\setminus\left(\bigcup_{U\in\cU}\mathcal{P}_U\right)\] be the set of global exceptions. 
    
    For each $d\in\mathcal{E}\setminus\{0\}$, apply Lemma \ref{lem:amend} once to obtain $\alpha_d\in\Q_{>0}$ and $U_d\in\cU$ with \[\frac1{\Prog(\alpha_d,1/d)}\subseteq\cS_k(U_d).\]
    For each $U\in\cU$, there are finitely many $d\in\mathcal{E}\setminus\{0\}$ such that $U_d=U$, because if this is so then $d\in\cS_k)$ but $d\notin\mathcal{P}_U$, and so $d\in \mathcal{E}_U$ and this set is finite. Therefore, for each $U\in\cU$, we define \[\Tilde{\mathcal{P}}_U=\mathcal{P}_U\cup\left(\bigcup_{d\in\mathcal{E}\setminus\{0\},U_d=U}\frac1{\Prog(\alpha_d,1/d)}\right)\] and it is still a finite union of reciprocals of arithmetic progressions. Then we have $\Tilde{\mathcal{P}}_U\subseteq\cS_k(U)$ and $\cS_k(\T^n)\setminus\{0\}=\bigcup_{U\in\cU}\Tilde{\mathcal{P}}_U$.
\end{proof}

\section*{Disclosure on the use of AI}
The idea of proving a generalization of \cite[Theorem 1.1]{jain_relative_2026} by looking at general rational piecewise-affine functions and applying the theory of linear programming came from the author. ChatGPT-5.6 and ChatGPT-6 were then used to carry out this strategy and provide the intricate proofs, which were simplified and rewritten by the author. 

\section*{Acknowledgements}
I am most thankful to Noah Kravitz, my advisor, for his guidance, suggestions, and encouragement throughout the time I was working on this problem; to say that he was impactful in my mathematical journey would be an understatement. I thank the organizers of the Lonely Runner Workshop at the University of Rostock for their hospitality and generosity. I am funded by the King's Scholarship (Thailand) and the EPSRC Vacation Internship scheme at the Mathematical Institute, University of Oxford.

\bibliographystyle{abbrv}
\bibliography{refs}
\end{document}